\documentclass[12pt]{amsart}

\usepackage{amsmath}
\usepackage{amssymb}
\usepackage{graphicx}
\usepackage{hyperref}
\usepackage{ifthen}
\usepackage[bbgreekl]{mathbbol}
\usepackage{pict2e}
\usepackage{xargs}
\usepackage{xspace}

\DeclareSymbolFontAlphabet{\amsmathbb}{AMSb}

\newcommand{\AD}{\mathtt{AD}}
\newcommand{\Bairespace}[1][]{
  \ifthenelse{\equal{#1}{}}{\functions{\N}{\N}}{\functions{#1}{\N}}
}

\newcommand{\bbs}{\mathbb{s}}
\newcommand{\calN}{\mathcal{N}}
\newcommand{\Cantorspace}[1][]{
  \ifthenelse{\equal{#1}{}}{\functions{\N}{2}}{\functions{#1}{2}}
}
\newcommand{\Cantortree}{\functions{<\N}{2}}
\newcommand{\composition}{\circ}
\newcommandx{\concatenation}[2][1 = undefined, 2 = undefined]{
  \ifthenelse{\equal{#1}{undefined}}{{}\smallfrown}{
    \ifthenelse{\equal{#2}{undefined}}{\bigoplus #1}{\bigoplus_{#1} #2}
  }
}
\newcommand{\DC}{\mathtt{DC}}
\newcommand{\definedterm}[1]{\textit{#1}}
\newcommandx{\Deltaclass}[2][1=,2=]{
  \ifthenelse{\equal{#2}{}}
  {\mathbf{\Delta}_{#1}}{\mathbf{\Delta}^{#1}_{#2}}
}
\newcommand{\equivalenceclass}[2]{[#1]_{#2}}
\newcommand{\equivalent}[1]{\equiv_{#1}}
\newcommand{\existnonmeagerlymany}{\exists^*}
\newcommand{\extendedby}{\sqsubseteq}
\newcommand{\extensions}[1]{\calN_{#1}}
\newcommand{\Ezero}{\amsmathbb{E}_0}
\newcommand{\forcomeagerlymany}{\forall^*}
\newcommand{\from}{\colon}
\newcommandx{\functions}[3][3 =]{
  \ifthenelse{\equal{#3}{}}{#2^{#1}}{#2_{#3}^{#1}}
}

\newcommand{\Gzero}{\amsmathbb{G}_0}
\newcommand{\image}[2]{#1(#2)}
\newcommand{\incomparable}[1]{\perp_{#1}}
\newcommandx{\intersection}[2][1 =, 2 =]{
  \ifthenelse{\equal{#1}{}}{\cap}{
    \ifthenelse{\equal{#2}{}}{\bigcap #1}{{\bigcap_{#1} #2}}
  }
}
\newcommand{\length}[1]{|#1|}
\newcommand{\mathand}{\text{ and }}
\newcommand{\N}{\amsmathbb{N}}
\newcommand{\pair}[2]{(#1, #2)}
\newcommand{\PD}{\mathtt{PD}}
\newcommandx{\Piclass}[2][1=,2=]{
  \ifthenelse{\equal{#2}{}}{\mathbf{\Pi}_{#1}}{\mathbf{\Pi}^{#1}_{#2}}
}

\newcommand{\predecessorordinal}[1]{\boldsymbol{\kappa}^1_{#1}}
\newcommand{\preimage}[2]{#1^{-1}(#2)}
\newcommandx{\product}[2][1 =, 2 =]{
  \ifthenelse{\equal{#1}{}}{\times}{
    \ifthenelse{\equal{#2}{}}{\prod #1}{{\prod_{#1} #2}}
  }
}
\newcommand{\projectiveordinal}[1]{\boldsymbol{\delta}^1_{#1}}
\renewcommand{\restriction}[2]{#1 \upharpoonright #2}

\newcommandx{\sequence}[2][2 = undefined]{
  \ifthenelse{\equal{#2}{undefined}}{(#1)}{
    (#1)_{#2}
  }
}
\newcommandx{\set}[2][2 = undefined]{
  \ifthenelse{\equal{#2}{undefined}}{\{ #1 \}}{
    \{ #1 \suchthat #2 \}
  }
}
\newcommand{\setcomplement}[1]{\mathord{\sim} #1}
\newcommandx{\Sigmaclass}[2][1=,2=]{
  \ifthenelse
    {\equal{#2}{}}
    {\mathbf{\Sigma}_{#1}}
    {\mathbf{\Sigma}^{#1}_{#2}}
}
\newcommand{\successor}[1]{#1^+}
\newcommand{\suchthat}{\mid}

\newcommandx{\union}[2][1 =, 2 =]{
  \ifthenelse{\equal{#1}{}}{\cup}{
    \ifthenelse{\equal{#2}{}}{\bigcup #1}{{\bigcup_{#1} #2}}
  }
}
\newcommand{\ZF}{\mathtt{ZF}}

\newcommand{\Baire}{Baire\xspace}
\newcommand{\Borel}{Bor\-el\xspace}
\newcommand{\Dilworth}{Dil\-worth\xspace}
\newcommand{\Harrington}{Har\-ring\-ton\xspace}
\newcommand{\Hausdorff}{Haus\-dorff\xspace}
\newcommand{\Kanovei}{Kan\-o\-vei\xspace}
\newcommand{\Kuratowski}{Kur\-at\-ow\-ski\xspace}
\newcommand{\Marker}{Mar\-ker\xspace}
\newcommand{\Mueller}{M\"{u}l\-ler\xspace}

\newcommand{\Shelah}{Shel\-ah\xspace}
\newcommand{\Souslin}{Sous\-lin\xspace}
\newcommand{\Ulam}{U\-lam\xspace}

\newenvironment{lemmaproof}{
  
  \begin{proof}
}{\end{proof}}

\newenvironment{theoremproof}{
  
  \begin{proof}
}{\end{proof}}

\newtheorem{lemma}{Lemma}
\newtheorem{theorem}[lemma]{Theorem}

\theoremstyle{definition}
\newtheorem*{acknowledgements}{Acknowledgements}

\title
  [Large antichains]
  {On the existence of large antichains for definable quasi-orders}

\author[B.D. Miller]{Benjamin D. Miller}
\address{
  Benjamin D. Miller \\
  Kurt G\"{o}del Research Center for Mathematical Logic \\
  Universit\"{a}t Wien \\
  W\"{a}hringer Stra{\ss}e 25 \\
  1090 Wien \\
  Austria
}
\email{benjamin.miller@univie.ac.at}
\urladdr{http://www.logic.univie.ac.at/benjamin.miller}

\author[Z. Vidny\'{a}nszky]{Zolt\'{a}n Vidny\'{a}nszky}
\address{
  Zolt\'{a}n Vidny\'{a}nszky \\
  Kurt G\"{o}del Research Center for Mathematical Logic \\
  Universit\"{a}t Wien \\
  W\"{a}hringer Stra{\ss}e 25 \\
  1090 Wien \\
  Austria
}
\email{zoltan.vidnyanszky@univie.ac.at}
\urladdr{
	http://www.logic.univie.ac.at/~vidnyanszz77/
}
	
\thanks{
  Both authors were supported in part by FWF Grants P28153 and
  P29999, and the second in part by the National Research,
  Development, and Innovation Office -- NKFIH Grants 113047, 104178,
  and 124749.
}
	
	

\begin{document}

\maketitle
	
\begin{abstract}
  We generalize
  \Harrington-\Marker-\Shelah's \Dilworth-style characterization of the
  existence of non-empty perfect antichains to co-analytic quasi-orders,
  establish the analogous theorem at the next definable cardinal, and
  consider generalizations beyond the first level of the projective
  hierarchy.
\end{abstract}

\section*{Introduction}

A \definedterm{quasi-order} is a reflexive transitive binary relation.
Associated with every such relation $R$ on a set $X$ are the
equivalence relation $x \mathrel{\equivalent{R}} y \iff (x \mathrel{R} y
\mathand y \mathrel{R} x)$ and the incomparability relation $x \mathrel
{\incomparable{R}} y \iff (\neg x \mathrel{R} y \mathand \neg y \mathrel
{R} x)$. We say that a set $Y \subseteq X$ is an \definedterm
{$R$-antichain} if $\restriction{R}{Y}$ is the diagonal on $Y$, and an
\definedterm{$R$-chain} if $\restriction{\mathord{\incomparable{R}}}{Y}$
is empty. A subset of a topological space is \definedterm{perfect} if it is
closed and has no isolated points, \definedterm{\Borel} if it is in the
$\sigma$-algebra generated by the open sets, \definedterm{analytic} if it
is a continuous image of a closed subset of $\Bairespace$, and
\definedterm{co-analytic} if its complement is analytic. In \S\ref{perfect},
we generalize \cite[Theorem 5.1]{HMS} from \Borel to co-analytic
quasi-orders:

\begin{theorem} \label{perfect:main}
  Suppose that $X$ is a \Hausdorff space and $R$ is a co-analytic
  quasi-order on $X$. Then exactly one of the following holds:
  \begin{enumerate}
    \item The space $X$ is a union of countably-many \Borel $R$-chains.
    \item There is a non-empty perfect $R$-antichain.
  \end{enumerate}
\end{theorem}

\noindent
Our proof uses only \Baire category arguments and the $\Gzero$
dichotomy \cite[Theorem 6.3]{KST}, which itself has a classical proof
\cite{Miller}. An interesting new wrinkle is that, while such arguments
typically utilize just one application of the $\Gzero$ dichotomy, ours
requires infinitely many.

A \definedterm{homomorphism} from a binary relation $R$ on $X$ to a
binary relation $S$ on $Y$ is a function $\phi \from X \to Y$ such that
$\image{(\phi \times \phi)}{R} \subseteq S$, and a \definedterm
{reduction} of $R$ to $S$ is a homomorphism from $R$ to $S$ that is
also a homomorphism from $\setcomplement{R}$ to $\setcomplement
{S}$. A \Borel equivalence relation $E$ on an analytic \Hausdorff space
$X$ is \definedterm{smooth} if there is a \Borel-measurable reduction
of $E$ to equality on $\Cantorspace$, and an analytic set $A \subseteq
X$ is \definedterm{$E$-smooth} if $\restriction{E}{A}$ is smooth. In \S\ref
{nonsmooth}, we establish the analog of \cite[Theorem 5.1]{HMS} at the
next \Borel cardinal:

\begin{theorem} \label{nonsmooth:main}
  Suppose that $X$ is an analytic \Hausdorff space and $R$ is a \Borel
  quasi-order on $X$. Then exactly one of the following holds:
  \begin{enumerate}
    \item There is a smooth \Borel superequivalence relation of
      $\equivalent{R}$ whose equivalence classes are $R$-chains.
    \item There is an $\equivalent{R}$-non-smooth perfect set whose
      quotient by $\equivalent{R}$ is an $(R / \mathord{\equivalent
      {R}})$-antichain.
  \end{enumerate}
\end{theorem}

\noindent
Our proof uses only \Baire category arguments and the $\Ezero$
dichotomy \cite[Theorem 1.1]{HKL}, which itself has a classical proof
\cite{Miller}, and reveals that the theorem holds for the rather degenerate
reason that its two alternatives are equivalent to those of the $\Ezero$
dichotomy (for $\equivalent{R}$).

A subset of a topological space $X$ is \definedterm{$\kappa$-\Borel} if it
is in the $\kappa$-complete algebra generated by the open sets,
\definedterm{$\kappa$-\Souslin} if it is a continuous image of a closed
subset of $\functions{\N}{\kappa}$, \definedterm{co-$\kappa$-\Souslin} if
its complement is $\kappa$-\Souslin, \definedterm{bi-$\kappa$-\Souslin}
if it is both $\kappa$-\Souslin and co-$\kappa$-\Souslin, and
\definedterm{$\aleph_0$-universally \Baire} if its pre-image under every
continuous function $\phi \from \Cantorspace \to X$ has the \Baire
property. Let $\Ezero$ denote the equivalence relation on
$\Cantorspace$ given by $c \mathrel{\Ezero} d \iff \exists n \in \N \forall
m \ge n \ c(m) = d(m)$. An \definedterm{embedding} is an injective
reduction. In \S\ref{generalizations}, we note that our arguments also
yield:

\begin{theorem} \label{generalizations:perfect:souslin}
  Suppose that $\kappa$ is an aleph, $X$ is a \Hausdorff space, and $R$
  is an $\aleph_0$-universally-\Baire co-$\kappa$-\Souslin quasi-order
  on $X$. Then at least one of the following holds:
  \begin{enumerate}
    \item The space $X$ is a union of $\kappa$-many $R$-chains.
    \item There is a non-empty perfect $R$-antichain.
  \end{enumerate}
\end{theorem}

\begin{theorem} \label{generalizations:nonsmooth:souslin}
  Suppose that $\kappa$ is an aleph, $X$ is a \Hausdorff space, and
  $R$ is an $\aleph_0$-universally-\Baire bi-$\kappa$-\Souslin
  quasi-order on $X$. Then at least one of the following holds:
  \begin{enumerate}
    \item There is a homomorphism from $\equivalent{R}$ to equality on
      $\functions{\kappa}{2}$ such that the pre-image of every singleton is
      an $R$-chain.
    \item There is a continuous embedding $\pi \from \Cantorspace \to X$
      of $\Ezero$ into $\equivalent{R}$ such that $\image{\pi}{\Cantorspace}
      / \mathord{\equivalent{R}}$ is an $(R / \mathord{\equivalent{R}})$-antichain.
  \end{enumerate}
\end{theorem}

\noindent
In the special case that $X$ is analytic, one can use the (non-classical)
arguments of \cite{Kanovei} to establish the strengthenings in which the
objects in condition (1) of these results are $\successor
{\kappa}$-\Borel/$\successor{\kappa}$-\Borel measurable.

A subset of an analytic \Hausdorff space is $\Sigmaclass[1][1]$ if it is
analytic, $\Piclass[1][n]$ if its complement is $\Sigmaclass[1][n]$,
$\Sigmaclass[1][n+1]$ if it is a continuous image of a $\Piclass[1][n]$
set, and $\Deltaclass[1][n]$ if it is both $\Piclass[1][n]$ and $\Sigmaclass
[1][n]$. Let $\projectiveordinal{n}$ denote the supremum of the lengths
of all $\Deltaclass[1][n]$ pre-wellorderings of $\Bairespace$. We say that
a $\Deltaclass[1][2n+1]$ equivalence relation $E$ on an analytic
\Hausdorff space $X$ is \definedterm{smooth} if there exists $\kappa <
\projectiveordinal{2n+1}$ for which there is a $\Deltaclass[1][2n+
1]$-measurable reduction of $E$ to equality on $\functions{\kappa}{2}$,
and an analytic set $A \subseteq X$ is \definedterm{$E$-smooth} if
$\restriction{E}{A}$ is smooth. Taking the known structure theory of the
projective sets as a black box, we note that our arguments also provide
classical proofs of the relevant special cases of the \Kanovei-style
strengthenings of Theorems \ref{generalizations:perfect:souslin} and \ref
{generalizations:nonsmooth:souslin} necessary to obtain:

\begin{theorem}[$\AD$] \label{generalizations:perfect:projective}
  Suppose that $n \in \N$, $X$ is an analytic \Hausdorff space, and $R$
  is a $\Piclass[1][2n+1]$ quasi-order on $X$. Then exactly one of the
  following holds:
  \begin{enumerate}
    \item The space $X$ is a union of $(<\projectiveordinal
      {2n+1})$-many $\Deltaclass[1][2n+1]$ $R$-chains.
    \item There is a non-empty perfect $R$-antichain.
  \end{enumerate}
\end{theorem}

\begin{theorem}[$\AD$] \label{generalizations:nonsmooth:projective}
  Suppose that $n \in \N$, $X$ is an analytic \Hausdorff space, and $R$
  is a $\Deltaclass[1][2n+1]$ quasi-order on $X$. Then exactly one of the
  following holds:
  \begin{enumerate}
    \item There is a smooth $\Deltaclass[1][2n+1]$ superequivalence
      relation of $\equivalent{R}$ whose equivalence classes are
      $R$-chains.
    \item There is an $\equivalent{R}$-non-smooth perfect set whose
      quotient by $\equivalent{R}$ is an $(R / \mathord{\equivalent
      {R}})$-antichain.
  \end{enumerate}
\end{theorem}

\noindent
In a future paper with \Mueller, we will establish the version of the
$\Gzero$ dichotomy necessary to obtain such results from $\PD$.

We work in $\ZF + \DC$ throughout.
 
\section{Perfect antichains} \label{perfect}

For each discrete set $D$ and sequence $s \in \functions{<\N}{D}$, we
use $\extensions{s}$ to denote the basic open set consisting of all
extensions of $s$ in $\functions{\N}{D}$. We use the notation
$\forcomeagerlymany x \in X \ P(x)$ to indicate that $\set{x \in X}[\neg
P(x)]$ is meager, and $\existnonmeagerlymany x \in X \ P(x)$ to indicate
that $\set{x \in X}[P(x)]$ is non-meager. Fix sequences $\bbs_n \in
\Cantorspace[n]$ such that $\forall s \in \Cantortree \exists n \in \N
\ s \extendedby \bbs_n$, and define $\Gzero = \union[n \in \N][{\set{\pair
{\bbs_n \concatenation \sequence{i} \concatenation c}{\bbs_n
\concatenation \sequence{1 - i} \concatenation c}}[c \in \Cantorspace
\mathand i < 2]}]$. While our proof of the characterization of the
existence of a non-empty perfect antichain requires infinitely-many
applications of the $\Gzero$ dichotomy, we need only one to establish
the following:

\begin{theorem} \label{perfect:weak}
  Suppose that $X$ is a \Hausdorff space, $R$ is a co-analytic
  quasi-order on $X$, and $X$ is not a union of countably-many \Borel
  $R$-chains. Then there are compact sets $K_i \subseteq X$ that are
  not unions of countably-many \Borel $R$-chains such that $\product
  [i < 2][K_i] \subseteq \mathord{\incomparable{R}}$.
\end{theorem}

\begin{theoremproof}
  As $\incomparable{R}$ is analytic and $X$ is not a union of
  countably-many \Borel $R$-chains, the $\Gzero$ dichotomy yields a
  continuous homomorphism $\phi \from \Cantorspace \to X$ from
  $\Gzero$ to $\incomparable{R}$. As the set $R_0 = \preimage{(\phi
  \times \phi)}{R}$ is co-analytic, it has the \Baire property (see, for
  example, \cite[21.6]{Kechris}), thus so too does $\incomparable{R_0}$.
  
  \begin{lemma} \label{perfect:weak:nonmeager}
    The relation $\incomparable{R_0}$ is non-meager.
  \end{lemma}
  
  \begin{lemmaproof}
    Suppose, towards a contradiction, that $\incomparable{R_0}$ is
    meager, and fix non-empty open sets $U_i \subseteq \Cantorspace$
    for which $R_0$ is comeager in $\product[i < 2][U_i]$ (see, for
    example, \cite[Proposition 8.26]{Kechris}). The \Kuratowski-\Ulam
    theorem (see, for example, \cite[Theorem 8.41]{Kechris}) ensures that
    the sets $C_0 = \set{c \in \Cantorspace}[\existnonmeagerlymany d \in
    U_0 \ c \mathrel{R_0} d]$ and $C_1 = \set{d \in \Cantorspace}
    [\forcomeagerlymany c \in U_0 \ c \mathrel{R_0} d]$ have comeager
    union, and \cite[Theorem 16.1]{Kechris} and the \Kuratowski-\Ulam
    theorem imply that they have the \Baire property. The
    \Kuratowski-\Ulam theorem also ensures that $C_0$ is non-meager,
    since otherwise $\forcomeagerlymany c, d \in U_0 \ (\neg c \mathrel
    {R_0} d \mathand d \mathrel{R_0} c)$, and $C_1$ is non-meager. As
    the $\Ezero$-saturation of every non-meager set with the \Baire
    property is comeager (see, for example, \cite[Theorem 8.47]{Kechris}),
    there are comeagerly-many $c \in \Cantorspace$ for which the sets
    $C_i \intersection \equivalenceclass{c}{\Ezero}$ non-trivially partition
    $\equivalenceclass{c}{\Ezero}$. As a straightforward induction reveals
    that $\Ezero$ is the equivalence relation generated by $\Gzero$, it
    follows that $(\product[i < 2][C_i]) \intersection \Gzero \neq \emptyset$,
    contradicting the fact that $\product[i < 2][C_i] \subseteq R_0$.  
  \end{lemmaproof}
  
  \begin{lemma} \label{perfect:weak:homomorphism}
    There are continuous homomorphisms $\phi_i \from \Cantorspace \to
    \Cantorspace$ from $\Gzero$ to itself for which $\product[i < 2][\image
    {\phi_i}{\Cantorspace}] \subseteq \mathord{\incomparable{R_0}}$.
  \end{lemma}
  
  \begin{lemmaproof}
    By Lemma \ref{perfect:weak:nonmeager}, there are non-empty open
    sets $U_i \subseteq \Cantorspace$ and dense open sets $V_n
    \subseteq \product[i < 2][U_i]$ such that $\intersection[n \in \N][V_n]
    \subseteq \mathord{\incomparable{R_0}}$. Recursively construct
    $u_{i,n} \in \Cantortree$ and $k_{i,n} \in \N$ such that $\product[i < 2]
    [\extensions{\phi_{i,n}(t_i)}] \subseteq V_n$ for all $t_0, t_1 \in
    \Cantorspace[n]$ and $\phi_{i,n}(\bbs_n) = \bbs_{k_{i,n}}$ for all $i <
    2$, where $\phi_{i,n} \from \Cantorspace[n] \to \Cantortree$ is given by
    $\phi_{i,n}(t) = u_{i,0} \concatenation \concatenation[m < n][\sequence
    {t(m)} \concatenation u_{i, m+1}]$. Then the functions $\phi_i \from
    \Cantorspace \to \Cantorspace$ given by $\phi_i(c) = \union[n \in \N]
    [\phi_{i,n}(\restriction{c}{n})]$ are as desired.
  \end{lemmaproof}
  
  It only remains to observe that if the functions $\phi_i$ are as in Lemma
  \ref{perfect:weak:homomorphism}, then the sets $K_i = \image{(\phi
  \composition \phi_i)}{\Cantorspace}$ are as desired.
\end{theoremproof}

We now establish our characterization of the existence of a non-empty
perfect antichain:

\begin{theoremproof}[Proof of Theorem \ref{perfect:main}]
  Conditions (1) and (2) are clearly mutually exclusive. To see $\neg(1)
  \implies (2)$, note that if condition (1) fails, then $X$ is the projection
  onto either coordinate of the complement of $\equivalent{R}$, thus
  analytic. Fix a continuous surjection $\phi \from \Bairespace \to X$, and
  recursively appeal to Theorem \ref{perfect:weak} to obtain functions
  $\psi_n \from \Cantorspace[n] \to \Bairespace[n]$ and sequences
  $\sequence{F_s}[{s \in \Cantorspace[n]}]$ of closed subsets of $X$ with
  the following properties:
  \begin{enumerate}
    \item $\forall s \in \Cantorspace[n] \ F_s$ is not a union of
      countably-many \Borel $R$-antichains.
    \item $\forall s \in \Cantorspace[n] \ F_s \subseteq \image{\phi}
      {\extensions{\psi_n(s)}}$.
    \item $\forall s \in \Cantorspace[n] \ F_{s \concatenation \sequence{0}}
      \union F_{s \concatenation \sequence{1}} \subseteq F_s$.
    \item $\forall s \in \Cantorspace[n] \ F_{s \concatenation
      \sequence{0}} \times F_{s \concatenation \sequence{1}} \subseteq
        \mathord{\incomparable{R}}$.
    \item $\forall i < 2 \forall s \in \Cantorspace[n] \ \psi_n(s) \extendedby
      \psi_{n+1}(s \concatenation \sequence{i})$.
  \end{enumerate}
  Define $\psi \from \Cantorspace \to \Bairespace$ by $\psi(c) = \union[n
  \in \N][\psi_n(\restriction{c}{n})]$, as well as $\pi = \phi \composition
  \psi$, noting that $\pi(c) \in \intersection[n \in \N][F_{\restriction{c}{n}}]$
  for all $c \in \Cantorspace$. To see that $\image{\pi}{\Cantorspace}$ is
  a perfect $R$-antichain, observe that if $c, d \in \Cantorspace$ are
  distinct, then there is a maximal natural number $n \in \N$ for which
  $\restriction{c}{n} = \restriction{d}{n}$, so the fact that $\pi(c) \in F_{s
  \concatenation \sequence{c(n)}}$ and $\pi(d) \in F_{s \concatenation
  \sequence{d(n)}}$ ensures that $\pi(c) \incomparable{R} \pi(d)$.
\end{theoremproof}

\section{Non-smooth antichains} \label{nonsmooth}

We now establish our characterization of the existence of a non-smooth
perfect set whose quotient is an antichain:

\begin{theoremproof}[Proof of Theorem \ref{nonsmooth:main}]
  To see that conditions (1) and (2) are mutually exclusive, note that if
  $E$ is a \Borel superequivalence relation of $\equivalent{R}$ whose
  classes are $R$-chains, and $A \subseteq X$ is a set whose quotient
  by $\mathord{\equivalent{R}}$ is an $(R / \mathord{\equivalent
  {R}})$-antichain, then $\restriction{\mathord{\equivalent{R}}}{A} =
  \restriction{E}{A}$. When $A$ is analytic, it follows that if $E$ is smooth,
  then so too is $\restriction{\mathord{\equivalent{R}}}{A}$.

  To see $\neg(1) \implies (2)$, note that if (1) fails, then $\equivalent{R}$
  is non-smooth, so the $\Ezero$ dichotomy yields a continuous
  embedding $\phi \from \Cantorspace \to X$ of $\Ezero$ into
  $\equivalent{R}$. As the set $R_0 = \preimage{(\phi \times \phi)}{R}$
  is \Borel, it has the \Baire property, thus so too does $\incomparable
  {R_0}$.
  
  \begin{lemma} \label{nonsmooth:main:incomparable}
    The relation $\incomparable{R_0}$ is comeager.
  \end{lemma}
  
  \begin{lemmaproof}
    If there exist $n \in \N$ and $s, t \in \Cantorspace[n]$ for which $R_0$
    is comeager in $\extensions{s} \times \extensions{t}$, then the fact
    that $\Ezero \subseteq R_0$ ensures that $R_0$ is comeager in
    $\extensions{s'} \times \extensions{t'}$ for all $s', t' \in \Cantorspace
    [n]$, and therefore comeager, thus so too is $\equivalent{R_0}$,
    contradicting the fact that the latter set is $\Ezero$.
  \end{lemmaproof}
  
  \begin{lemma} \label{nonsmooth:main:homomorphism}
    There is a continuous embedding $\psi \from \Cantorspace \to
    \Cantorspace$ of $\Ezero$ into itself that is also a homomorphism
    from $\setcomplement{\Ezero}$ to $\incomparable{R_0}$.
  \end{lemma}
  
  \begin{lemmaproof}
    By Lemma \ref{nonsmooth:main:incomparable}, there are dense open
    sets $U_n \subseteq \Cantorspace \times \Cantorspace$ such that
    $\intersection[n \in \N][U_n] \subseteq \mathord{\incomparable{R_0}}$.
    We can clearly assume that these sets are decreasing and disjoint
    from the diagonal. Recursively construct $u_{i,n} \in \Cantortree$ such
    that $\length{u_{0,n}} = \length{u_{1,n}}$ and $\product[i < 2]
    [\extensions{\psi_{n+1}(t_i \concatenation \sequence{i})}] \subseteq
    U_n$ for all $t_0, t_1 \in \Cantorspace[n]$, where $\psi_{n+1} \from
    \Cantorspace[n+1] \to \Cantortree$ is given by $\psi_{n+1}(t) =
    \concatenation[m \le n][\sequence{t(m)} \concatenation u_{t(m),m}]$.
    Then the map $\psi \from \Cantorspace \to \Cantorspace$ given by
    $\psi(c) = \union[n \in \N][\psi_n(\restriction{c}{n})]$ is as desired.
  \end{lemmaproof}
  
  It only remains to observe that if the function $\psi$ is as in Lemma
  \ref{nonsmooth:main:homomorphism}, then the set $\image
  {(\phi \composition \psi)}{\Cantorspace}$ is as desired.
\end{theoremproof}

\section{Generalizations} \label{generalizations}

Simplifications of the classical proofs of the $\Gzero$ and $\Ezero$
dichotomies can be used to obtain a continuous homomorphism from
$\Gzero$ to every $\kappa$-\Souslin graph on a \Hausdorff space with
no $\kappa$-coloring, and a continuous embedding of $\Ezero$ into
every bi-$\kappa$-\Souslin equivalence relation on a \Hausdorff space
which is not reducible to equality on $\functions{\kappa}{2}$ \cite{Miller}.
By using these facts in lieu of the usual dichotomies in our proofs of
Theorems \ref{perfect:main} and \ref {nonsmooth:main}, we obtain
proofs of Theorems \ref{generalizations:perfect:souslin} and \ref
{generalizations:nonsmooth:souslin}.

If $\AD$ holds and $n \in \N$, then a subset of an analytic \Hausdorff
space is $\Deltaclass[1][2n+1]$ if and only if it is $\projectiveordinal{2n+
1}$-\Borel \cite{Martin, Moschovakis}, and there is a cardinal
$\predecessorordinal{2n+1}$ for which $\projectiveordinal{2n+1} = 
\successor{(\predecessorordinal{2n+1})}$ \cite{Kechris:projective} and a
subset of an analytic \Hausdorff space is $\Sigmaclass[1][2n+1]$ if and
only if it is $\predecessorordinal{2n+1}$-\Souslin (see, for example, \cite
[Theorem 2.21]{Jackson}). It follows that continuous images of
$\projectiveordinal{2n+1}$-\Borel sets are $\predecessorordinal{2n+
1}$-\Souslin, a fact which alone ensures that the classical proofs of the
$\Gzero$ and $\Ezero$ dichotomies yield the special cases of the
\Kanovei-style generalizations thereof at $\predecessorordinal{2n+1}$.
By using these in lieu of the usual dichotomies in our proofs of
Theorems \ref{perfect:main} and \ref{nonsmooth:main}, we obtain
proofs of the \Kanovei-style strengthenings of Theorems \ref
{generalizations:perfect:souslin} and \ref
{generalizations:nonsmooth:souslin} at $\predecessorordinal{2n+1}$.
As $\AD$ also ensures that every subset of a topological space is
$\aleph_0$-universally \Baire (see, for example, \cite[Theorem 38.17]
{Kechris}), Theorems \ref{generalizations:perfect:projective} and \ref
{generalizations:nonsmooth:projective} easily follow.

\begin{acknowledgements}
  We would like to thank Sandra \Mueller for a useful conversation
  concerning the results of \S\ref{generalizations}.
\end{acknowledgements}

\bibliographystyle{amsalpha}
\bibliography{hms}

\end{document}